\PassOptionsToPackage{unicode}{hyperref}
\PassOptionsToPackage{naturalnames}{hyperref}

\documentclass[a4paper,11pt]{article}

\usepackage{amsfonts,amssymb,amsmath,amsthm}	
\usepackage[english]{babel}                 	
\usepackage[T1]{fontenc}                    	
\usepackage{setspace}      						
\usepackage{graphicx}      						
\usepackage{tikz}          						
\usepackage{enumitem}
\usepackage[permil]{overpic}       				
\usepackage{fancybox}							
\usepackage{ifthen}								
\usepackage{subcaption}
\usepackage[ plainpages = false, pdfpagelabels,
	     bookmarks,
	     bookmarksopen = true,
	     bookmarksnumbered = true,
	     breaklinks = true,
	     linktocpage,
	     pagebackref,         
	     colorlinks = true,  
         hypertexnames=false,   
	     linkcolor = green!50!black,
	     urlcolor  = blue!50!black,
	     citecolor = blue!50!black,
	     anchorcolor = blue,
	     hyperindex = true,
	     hyperfigures
	     ]{hyperref}
\usepackage{tikz}
\usepackage{tikz-3dplot}
\usepackage[numbers]{natbib}
\usepackage{bibentry} 
\usepackage{fancyhdr, graphicx}
\usepackage{pgfgantt}
\usepackage{bm}

\usepackage{bbm}
\usepackage[colorinlistoftodos,prependcaption]{todonotes}
\usepackage[capitalise]{cleveref}
\usepackage{autonum} 
\usepackage{amssymb}
\usepackage{titlesec}
\usepackage{booktabs}
\usepackage[parfill]{parskip} 
\usepackage{algorithm}
\usepackage{algorithmic}
\usepackage{pgfplots}
\pgfplotsset{compat=newest}
\usetikzlibrary{plotmarks}
\usetikzlibrary{arrows.meta}
\usepgfplotslibrary{patchplots}
\usepackage{grffile}
\usepackage{amsmath}
\pgfplotsset{plot coordinates/math parser=false}
\newlength\figureheight
\newlength\figurewidth
\pgfplotsset{every axis/.append style={
                    label style={font=\small},
                    tick label style={font=\footnotesize},
                    legend style={font=\tiny},
                    title style={font=\small}
                    }}
\usetikzlibrary{plotmarks}
\pgfplotsset{minor grid style={dotted,gray}} 
\pgfplotsset{every axis/.append style={thick, tick style=semithick}}
\pgfplotsset{every mark/.append style={mark size=1pt}}
\usepgfplotslibrary{groupplots}
\pgfplotsset{xticklabel style={/pgf/number format/fixed,
        /pgf/number format/precision=5},yticklabel style={/pgf/number format/fixed,
        /pgf/number format/precision=5}}

 \pgfplotsset{
compat=1.11,
legend image code/.code={
\draw[mark repeat=2,mark phase=2]
plot coordinates {
(0cm,0cm)
(0.15cm,0cm)        
(0.3cm,0cm)         
};%
}
}

\graphicspath{{./01_article/figures/}{./figures/}}
\newcommand{\includetikz}[1]{%
        \includegraphics{#1.pdf} 
}

\renewcommand{\d}{\, \mathrm{d}}                                   
\newcommand{\vect}[1]{\textrm{\boldmath${#1}$}}                 
\newcommand{\norm}[1]{\Vert{#1}\Vert}

\newcommand{\energy}{{E}}

\newcommand{\Etot}{\energy_\mathrm{tot}}                 
\newcommand{\Epot}{\energy_\mathrm{pot}}                 
\newcommand{\Ekin}{{\energy}_\mathrm{kin}}               
\newcommand{\Bfield}{\vect{B}}

\newcommand{\map}{\vect{X}}

\newcommand{\vmap}{\map}

\newcommand{\email}[1]{\href{mailto:#1}{\texttt{#1}}}           
\newcommand{\RomanNumeralCaps}[1]                               
    {\MakeUppercase{\romannumeral #1}}

\usepackage{bm}
\newcommand{\vu}{{\bm u}}
\newcommand{\vv}{{\bm v}}

\newcommand{\vtx}{{\omega}}
\newcommand{\vx}{{\bm x}}

\newcommand{\vB}{{\bm B}}
\newcommand{\vW}{{\bm W}}

\newcommand{\dt}{{\Delta t}}

\newcommand{\vX}{{\bm X}}

\newcommand{\vtX}{\tilde{\bm X}}
\newcommand{\vtu}{\tilde{\bm u}}

\newcommand{\vhX}{\mathcal{\bm X}}
\newcommand{\vthX}{\tilde{\mathcal{\bm X}}}
\newcommand{\F}{F}
\newcommand{\hF}{\mathcal{F}}
\newcommand{\tf}{\tilde{f}}

\newcommand{\adj}{\text{adj}}

\newcommand{\curl}{\nabla \times}
\renewcommand{\div}{\nabla \cdot}
\newcommand{\grad}{\nabla}
\newcommand{\Laplace}{\Delta}

\newcommand{\Ih}{\mathcal{I}}

\newcommand{\Lie}{\mathfrak{L}}

{\theoremstyle{definition}}

{\theoremstyle{remark} }            

{\theoremstyle{lemma}
}

{\theoremstyle{proposition}
\newtheorem{proposition}{Proposition}}

\graphicspath{{figures/}}

\binoppenalty=\maxdimen 
\relpenalty=\maxdimen   

\usepackage[margin=1in]{geometry}

\title{\textbf{A Characteristic Mapping Method with Source Terms:}\\
{\Large Applications to Ideal Magnetohydrodynamics}}
\author{
Xi-Yuan Yin \footnote{CNRS, Ecole Centrale de Lyon, INSA Lyon, Université Claude Bernard Lyon 1, LMFA, UMR5509, 69130, Écully, France}
\footnote{Current address: Max Planck Institute for Mathematics in the Sciences, 04103 Leipzig, Germany, \email{xi.yin@mis.mpg.de}},
Philipp Krah\footnote{Institut de Mathématiques de Marseille (I2M), Aix-Marseille Université, \email{philipp.krah@univ-amu.fr}},
Jean-Christophe Nave\footnote{McGill University, Montreal \email{jcnave@math.mcgill.ca}}\, and
Kai Schneider\footnote{Institut de Mathématiques de Marseille (I2M), Aix-Marseille Université, \email{kai.schneider@univ-amu.fr}}
}

\nobibliography{references} 
\begin{document}

\maketitle

\begin{abstract}

This work introduces a generalized characteristic mapping method designed to handle non-linear advection with source terms. The semi-Lagrangian approach advances the flow map, incorporating the source term via the Duhamel integral. We derive a recursive formula for the time decomposition of the map and the source term integral, enhancing computational efficiency. Benchmark computations are presented for a test case with an exact solution and for two-dimensional ideal incompressible magnetohydrodynamics (MHD). Results demonstrate third-order accuracy in both space and time. The submap decomposition method achieves exceptionally high resolution, as illustrated by zooming into fine-scale current sheets. An error estimate is performed and suggests third order convergence in space and time.

\end{abstract}


\section{Introduction}

Conservation laws with source terms have a vast field of applications, ranging from geophysical flows, via chemical reacting flows, collision terms in kinetic equations or electromagnetic forces in Maxwell equations. Classical numerical schemes suffer from the stiffness introduced by the source terms and special care must be taken to avoid the convergence properties being degraded or even lost. For further details, we refer to the textbooks of Toro~\cite{toro2013riemann} and Leveque~\cite{leveque1992numerical} and the cited literature.
Prominent examples of methods handling conservation laws with source terms are splitting schemes, or e.g. the finite volume evolution Galerkin method by Morton and co-workers, see e.g. \cite{lukacova2004finite}.

This paper aims to develop a numerical method for nonlinear advection problems with source terms. The underlying transport velocity is assumed to be divergence-free and the equations describe trajectories on the space of volume-preserving diffeomorphisms. While our approach is general, the present paper focuses on electrically conducting fluid flows, i.e. on ideal magnetohydrodynamics. 
Magnetohydrodynamics (MHD) describes the interaction between conducting fluids and magnetic fields. In the ideal setting, both the resistivity and the viscosity of the fluid vanish, which is called ideal MHD, see e.g. the textbook \cite{davidson2016introduction}.
The governing equations are the incompressible Euler equations coupled with Maxwell equations involving the Lorenz force.
Two-dimensional (2d) magnetohydrodynamics (MHD) shares a crucial feature with three-dimensional (3d) fluid turbulence: a direct energy cascade toward smaller scales, potentially leading to anomalous dissipation and the formation of singularities. The development of fine-scale or even singular structures—such as thin current sheets and the exponential growth of current density and vorticity in 2d incompressible MHD \cite{grauer1998geometry}—poses significant challenges for numerical methods.

Over the last decades numerous MHD solvers in different domains, for instance astrophysics, magnetically confined fusion, liquid metals, have been developed. A short overview can be found in the introduction of \cite{morales2014simulation}.
Fyfe, Montgomery \& Joyce~\cite{fyfe1977dissipative}, Pouquet~\cite{pouquet1978two} and Orszag-Tang \cite{orszag1979small}
proposed Fourier pseudo-spectral methods for solving the viscous and resistive 2d MHD equations in periodic geometry.
Small-scale structures in three-dimensional magnetohydrodynamic turbulence using a pseudospectral method have been studied in~\cite{mininni2006small}.
Schneider et al.~\cite{schneider2011pseudo} proposed the coupling of the pseudo-spectral method with volume penalization to compute MHD turbulence in confined domains.
Simulations of ideal 2d MHD using Fourier Galerkin truncated discretizations and analyzing thermalization following the ideas of TD Lee, i.e., energy equipartition of spectral approximations \cite{lee1952some}, have been carried out by Brachet's group~\cite{krstulovic2011alfven}.
Some review of MHD solvers developed to compute
fusion-plasma-related flows is given in \cite{jardin2012review}.
Furthermore, adaptive 2d and 3d MHD computations using multiresolution methods can be found in \cite{gomes2021adaptive} and computations using adaptive mesh refinement have been proposed by Grauer's group  \cite{friedel1997adaptive}. 

Recently, the characteristic mapping methods (CMM), a novel numerical framework for advection problems, have been developed for various homogeneous transport problems. This approach is based on computing the inverse flow map or back-to-label map generated by the transporting velocity field, from which transported quantities can be directly evaluated by pullback. So far, these methods have shown success in the discretization of homogeneous equations such as linear transport \cite{mercier2020characteristic,taylor2023projection}, the incompressible Euler equations in 2d periodic domain \cite{yin2021characteristic, bergmann2024singularity}, on the 2-sphere \cite{taylor2023characteristic}, in 3d periodic domain \cite{yin2023characteristic}, as well as the 1+1d Vlassov-Poisson equations \cite{KrahYinBergmannSchneiderNave2023}. These methods preserve the relabeling symmetry, resulting in non-diffusive transport: the modified equations are closer to the Langrangian-averaged Euler-$\alpha$ equations or the Kelvin-filtered models where regularization is only applied to the transporting velocity field. The group property of the diffeomorphisms also allows for a submap decomposition of a long-time flow resulting in efficient numerical resolution of exponential growth in scale separation. However, so far these schemes have been limited to homogeneous equations without source terms.

In this paper, we formulate the CMM for advection equations with source terms and focus on the extension of the method to the ideal magnetohydrodynamics (MHD) equations. The main difficulty in this extension is the presence of the Lorentz force as source term in the momentum equation and therefore requires the development of the CMM for inhomogeneous advection. Our approach is the following: by introducing the characteristic map, the nonlinear transport in the Euler and ideal MHD equations is rewritten in a ``quasilinear'' form. When going from the homogeneous to inhomogeneous equations, this allows for the application of Duhamel's principle where the characteristic maps serve as solution operator to the linearized equation (see also \cite{TaylorYinNave2024functional} for analytic considerations for this approach and \cite{yin2021thesis} for applications to other equations). The numerical method then consists in the evolution of two quantities: the backward characteristic map, using the same method as in the homogeneous case, and the Duhamel time integral of the source term, which we also call ``accumulated source term'', solved using a semi-Lagrangian method on a fixed Eulerian grid. Again by Duhamel's principle, a submap decomposition method for the CMM also applies to the source term time integral, and as long as remapping is performed before the accumulated source term saturates the numerical spatial resolution, it will not incur numerical artificial diffusion or thermalization errors. This approach results in a modification of the original CMM framework for general advection equations with source terms. In the case of the ideal MHD equations, the method is further simplified by the structure of the source term as the Lorentz force depends only on the magnetic field, which is a differential 2-form Lie-advected by the flow. This means that the magnetic field at any time instance can be directly expressed in terms of the initial magnetic field and the characteristic map, more precisely by differential pullback. This type of source term structures are known as advected parameters since the initial conditions become essentially ``parameters'' to the map evolution equation. Extensive studies of these equations in the framework of geometric hydrodynamics have been carried out in \cite{holm1998euler, khesin2021geometric, arnold2008topological}. Numerically, this implies that the inhomogeneous advection equation for the accumulated source term is also quasi-linear, thus avoiding difficulties with the stability of time integration schemes. For more general cases with fully nonlinear source term equations, we expect that modifications to the time evolution schemes are required.

The outline of the manuscript is as follows.
In Sec.~\ref{sec:CMM} we present the CMM and its extension for handling source terms.
Sec.~\ref{sec:MHD} recalls the governing equations of ideal MHD in vorticity formulation. The numerical implementation of the CMM is then exposed in Sec.~\ref{sec:numimp}. Numerical validation for different test cases is given in
Sec.~\ref{sec:numtest}.
Finally, some conclusions and future directions are discussed in Sec.~\ref{sec:concl}.

\section{Characteristic Mapping Method}
\label{sec:CMM}

In this section we recall the Characteristic Mapping Method (CMM) studied in \cite{yin2021characteristic,yin2023characteristic,bergmann2024singularity,KrahYinBergmannSchneiderNave2023} and describe an extension of the method to account for source terms, particularly for the ideal magnetohydrodynamics equations.

The CMM was successfully employed for the linear advection, incompressible Euler equations in two and three dimensions as well as on the sphere. The main idea of the approach is to compute the inverse flow map generated by the velocity field. This yields at all times a diffeomorphism on the domain between time the current and the initial times. For our current investigation, we will only use the 2d periodic square domain $\mathbb{T}^2$.

For a given time dependent velocity field $\vu : \mathbb{T}^2 \times \mathbb{R} \to \mathbb{R}^2$, we denote by $\vX_{[s, t]}$ the flow map of $\vu$ from time $s$ to $t$. Under suitable assumptions on $\vu$, the flow maps are diffeomorphism with the properties
\begin{equation} \label{eqn:mapProperties}
        \vX_{[t, t]} (\vx) = \vx, \qquad \vX_{[r, t]} \circ \vX_{[s, r]} = \vX_{[s, t]}, \qquad \vX_{[s, t]} = \vX_{[t, s]}^{-1},
\end{equation}
for $s, r, t$ arbitrary and satisfy the equations
\begin{subequations}
    \begin{align}
        & \partial_t \vX_{[s, t]} = \vu ( \vX_{[s, t]} ), \\
        & (\partial_t + \vu \cdot \grad) \vX_{[t, s]} = 0 .
    \end{align}
\end{subequations}

We note that the backward map $\vX_{[t, 0]}$ is the solution operator for the advection equation under the velocity field $\vu$ in the sense that for any scalar field $\theta_0$ satisfying
\begin{equation} \label{eqn:scalarAdv}
    (\partial_t + \vu \cdot \grad) \theta = 0, \qquad \theta(\vx, 0) = \theta_0,
\end{equation}
we have a direct representation of the solution at time $t$ by pullback:
\begin{equation} \label{eqn:scalarPullback}
    \theta (\cdot, t) = \vX_{[t, 0]}^* \theta_0 =  \theta_0 \circ \vX_{[t, 0]}  ,
\end{equation}
where the superscript asterisk denotes the pullback operator, in general the dual operator of a change of coordinates by the map. This property was used in \cite{yin2021characteristic} to solve the incompressible Euler equations in 2d using the vorticity formulation since vorticity, in the 2d case, is an advected scalar field.

The compositional properties \eqref{eqn:mapProperties} of the maps naturally lead to a geometric numerical scheme for advection problems. The general discretization strategy is described as follows. Let $\Ih_M$ be an interpolation operator associated with a given grid $M$ on the domain $\Omega$. We denote by $\vhX_{[t, s]}$ the numerical map, the time-stepping scheme can be summarized as
\begin{equation} \label{eq:CMM_timestep}
    \vhX_{[t+\dt, s ]} = \Ih_M \left[ \vhX_{[t, s]} \circ \vtX_{[t+\dt, t]}  \right] , \qquad \vhX_{[s, s]} = \textit{id} ,
\end{equation}
where the one-step map $\vtX_{[t+\dt, t]}$ is obtained from a $\dt$ backward in time integration of the velocity field $\vu$ applied on the ODE
\begin{equation}
    \partial_r \vtX_{[t+\dt, r ]} = -\vu (\vtX_{[t+\dt, r ]}, r ).
\end{equation}
The equation is in some sense linearized by computing the characteristic map. The nonlinearity is limited to the evaluation of the velocity field $\vu$ which depends on the initial conditions and the map $\vX_{[t, 0]}$. The exact form of this velocity depends on the equation studied and the presence of source terms, we will present this in the next section.

For a short time $t - \tau_{i-1}$, an appropriate grid $M$ will be sufficient to accurately approximate the exact map $\vX_{[t, \tau_{i-1}]}$. When the numerical resolution of the grid is exhausted, a remapping is triggered with $\tau_i = t$ as remapping time. We thereby employ a submap decomposition method to approximate the long-time map $\vX_{[t, 0]}$ by
\begin{equation}
    \vhX_{[t, 0]} = \vhX_{[\tau_1, 0]} \circ \vhX_{[\tau_2, \tau_1]} \circ \vhX_{[\tau_3, \tau_2]} \circ \dots \circ \vhX_{[\tau_{n-2}, \tau_{n-1}]} \circ \vhX_{[t, \tau_{n-1}]},
\end{equation}
for a subdivision $0 = \tau_0 < \tau_1 < \tau_2 < \ldots < \tau_{n-1} < \tau_n = t$ of the full $[0, t]$ time interval. Each submap $\vhX_{[\tau_i, \tau_{i-1}]}$ is independently computed using the scheme eq. \eqref{eq:CMM_timestep} and the remapping times $t_i$ are triggered dynamically as the resolution capacity of the grid is reached. We see that as long as the numerical resolution is sufficient, each numerical submap $\vhX_{[\tau_i, \tau_{i-1}]}$ remains a diffeomorphism and therefore, the pullback is invertible. Practically, this means that the solution operator $\vhX_{[t, 0]}^*$ achieves an arbitrary resolution of all advected quantities and the pullback remains an infinite dimensional operator despite being discretized on a finite-dimensional function space.

\subsection{Source Terms}

Consider the scalar advection equation \eqref{eqn:scalarAdv} with solution given in terms of the initial condition and the characteristic map \eqref{eqn:scalarPullback}. We note that the evolution of the scalar $\theta(t)$ is given by a linear action of the map in the sense that $\theta (t) = \vX_{[t, 0]}^* \theta_0 = \theta_0 \circ \vX_{[t, 0]}$ and, for $a,b$ arbitrary, $a \theta (t) + b \phi(t) = \vX_{[t, 0]}^* ( a \theta_0 + b \phi_0)$ for $\theta, \phi$ solutions to linear homogeneous advection equations with initial conditions $\theta_0, \phi_0$, under the same velocity field $\vu$. Solutions to the inhomogeneous equation can therefore be obtained from the homogeneous equations through Duhamel's principle. This in fact also holds for the Lie-advection of differential $k$-forms where $\vX_{[t, 0]}^*$ denotes the pullback operator, and more generally, for linear (right) actions of the diffeomorphism group. 

Therefore, for the inhomogeneous advection equation
\begin{equation} \label{eq:inhomogAdv}
    (\partial_t + \vu \cdot \nabla) \theta = f, \qquad \theta(\vx, 0) = \theta_0(\vx) ,
\end{equation}
we can perform the decomposition
\begin{equation}
    \theta(\cdot, t) = \vX_{[t, \tau_1]}^* \left( \vX_{[\tau_1, 0]}^* \theta_0 + \int_0^{\tau_1} \vX_{[\tau_1, s]}^* f(\cdot, s) d s  \right)  + \int_{\tau_1}^t \vX_{[t, s]}^* f(\cdot, s) d s ,
\end{equation}
by commuting pullback with time integration. Defining
\begin{equation} \label{eq:sourceAccumDef}
    \F_{[\tau_i, \tau_{i-1}]} = \int_{\tau_{i-1}}^{\tau_i} \vX_{[\tau_i, s]}^* f(\cdot, s) d s,
\end{equation}
we can write a recursive formula for the submap decomposition expression
\begin{equation}
    \theta(\cdot, \tau_i) = \vX_{[\tau_i, \tau_{i-1}]}^* \theta(\cdot, \tau_{i-1}) + \F_{[\tau_{i}, \tau_{i-1}]} ,
\end{equation}
for $0 = \tau_0 < \tau_1 < \tau_2 < \ldots < \tau_{n-1} < \tau_n = t$ and $\theta(\cdot, 0) = \theta_0 (\cdot)$. Here, we note that within one remapping subinterval, for $t \in [\tau_{i-1}, \tau_{i}]$, the source term accumulation $\F_{[t, \tau_{i-1}]}$ as a function of $t$ satisfies
\begin{equation} \label{eq:sourceTermAccumEq}
    (\partial_t + \vu  \cdot \nabla) \F_{[t, \tau_{i-1}]} = f(\cdot, t) , \qquad \F_{[\tau_{i-1}, \tau_{i-1}]} = 0 .
\end{equation}
Therefore, in the CMM for inhomogeneous advection we must at the same time discretize two quantities, the current submap $\vX_{[t, \tau_{i-1}]}$ and the current source term accumulation $F_{[t, \tau_{i-1}]}$. 

\section{Ideal Magnetohydrodynamics in Vorticity Formulation}
\label{sec:MHD}

One of the main goals of this paper is to extend the CMM to the ideal magnetohydrodynamic (MHD) equations in two-dimensional space, i.e. with vanishing resistivity and viscosity. Our numerical method uses the vorticity formulation for the momentum transport coupled with the magnetic field advection equation:
    \begin{align}
         &\partial_t \vtx + (\vu \cdot \grad) \vtx = \curl \left(\curl  \vB \times \vB \right), \\
         &\partial_t \vB + \curl (\vB \times \vu) = 0,\label{eqns:MHD}\\
         &\vu= -\curl \Laplace^{-1} \vtx, 
    \end{align}
where $\vu$ is the velocity field, $\vtx$ is the vorticity field, and $\vB$ is the magnetic field parallel to the plane. Numerical simulation of these equations is already challenging in 2d \cite{grauer1998geometry}. The Lie-advection of the magnetic field may result in its local intensification in a process similar to the vortex stretching effect in the 3d Euler equations. As a result, the solution may contain extremely small scales in the magnetic field observed as an exponential growth in the maximum current density \cite{grauer1998geometry}. Resolving these fine scales makes the numerical simulation of this system difficult and costly. However, owing to the compositional structure of the CMM, an efficient resolution of exponential growth in scales of the numerical solution is possible as observed in \cite{yin2021characteristic, yin2023characteristic}.

In order to apply the CMM to the ideal MHD equations, we must solve an additional transport equation for the magnetic field $\vB$, this equation can be written as a Lie-advection equation. We consider the 2d transport equations as a system in three-dimensional space where all components of the solution are invariant in the third direction $\bm{e}_3$. In this case, both the vorticity field and the magnetic fields are Lie-advected differential 2-forms. In terms of the corresponding vector fields, the $\vtx$ is perpendicular to the plane and $\vB$ is parallel to the plane. For a 2-form $\vW$, the Lie derivative with respect to a velocity field $\vu$ is given by Cartan's homotopy formula:
\begin{equation}
\Lie_\vu \vW = \curl ( \vW \times \vu) +  \vu (\div \vW ) = (\vu \cdot \grad) \vW - (\vW \cdot \grad )\vu + \vW (\div \vu ).
\end{equation}
We see that plugging $\vW = \vtx \bm{e}_3$ gives simply the scalar $\vu$ directional derivative, whereas for $\vW = \vB$, we obtain the induction term in the magnetic field equation.

Similar to the scalar advection case, the solution operator for the Lie-advection equation is the pullback operator by the backward characteristic map:
\begin{equation}
    \vW(\vx, t) = \vX_{[t, 0]}^* \vW_0 (\vx) = \adj (D \vX_{[t, 0]} ) ( \vW_0 \circ \vX_{[t, 0]} ) .
\end{equation}
where $\adj$ denotes the adjugate matrix. This pullback formula was used in \cite{yin2023characteristic} for the incompressible Euler equations in 3D where the vorticity is fully three-dimensional and the adjugate of the Jacobian matrix was responsible for the vortex stretching.

Here, in the 2d ideal MHD case, we note that when $\vW = \vtx \bm{e}_3$, $\vW_0$ is perpendicular to the plane and the adjugate Jacobian term drops out. However, for $\vW = \vB$, we obtain the evaluation formula for the magnetic field based on initial data and the backward map:
\begin{equation} \label{eqn:2formPullback}
    \vB ( \cdot, t) = \adj (D \vX_{[t, 0]} ) ( \vB_0 \circ \vX_{[t, 0]} ) .
\end{equation}

Writing the magnetic fields through the pullback operator greatly simplifies the numerical scheme as $\vB(\vx, t)$ can be directly expressed as a function of $\vB_0$ and the map $\vX_{[t, 0]}$ implying that no extra equations need to be solved to compute the magnetic field. The Lorentz force may then be thought of as a function of the backward map depending on $\vB_0$ which can be viewed as a parameter of the equations. This allows us to rewrite \eqref{eqns:MHD} in terms of the characteristic map:
\begin{subequations} \label{eqns:MHD_CMM}
    \begin{align}
        & \partial_t \vX_{[t, 0]} + (\vu \cdot \grad) \vX_{[t, 0]} = 0, \qquad \vX_{[0, 0]} (\vx) = \vx, \label{subeq:CMM} \\
        & \partial_t F_{[t,0]} + (\vu \cdot \grad) F_{[t, 0]} =  \curl \left(\curl  \vB \times \vB \right), \qquad F_{[0,0]} = 0 , \label{subeq:sourceAccum} \\
        & \vtx (\cdot, t) =  \vtx_0 \circ \vX_{[t, 0]} + F_{[t,0]} , \label{subeq:vortEval} \\
        & \vB (\cdot, t) = \adj (D \vX_{[t,0]}) \vB_0 \circ \vX_{[t, 0]} , \label{subeq:magneticField_pullback} \\
        & \vu = - \curl \Laplace^{-1} \vtx .
    \end{align}
\end{subequations}

We note that equations \eqref{subeq:CMM} and \eqref{subeq:sourceAccum} are inviscid transport equations written in the Eulerian frame. There is no \emph{a priori} bound on the small scales generated by these equations and hence Eulerian discretizations are only valid for limited time intervals until the grid resolution is no longer sufficient. To overcome this problem, the CMM allows for a submap decomposition method where a long-time map $\vX_{[t, 0]}$ can be represented numerically as a composition of many submaps over shorter times wherein each submap has sufficient spatial resolution. Let $0 = \tau_0 < \tau_1 < \ldots < \tau_{k-1} < \tau_k = t$ be a subdivision of the time interval $[0, t]$, the following recursion gives an evaluation method for $\vX_{[t, 0]}$ and $F_{[t, 0]}$ amenable for numerical implementation:
\begin{subequations}
    \begin{align}
        & \vX_{[t, \tau_i]} = \vX_{[\tau_{i+1}, \tau_i]} \circ \vX_{[t, \tau_{i+1}]}, \label{subeq:mapRecursion} \\
        & F_{[t, \tau_i]} = F_{[t, \tau_{i+1}]} + F_{[\tau_{i+1}, \tau_i]} \circ \vX_{[t, \tau_{i+1}]} \label{subeq:sourceRecursion} .
    \end{align}
\end{subequations}

As a result, with this decomposition of the time interval, only the map and source term $\vX_{[t, \tau_i]}$ and $F_{[t, \tau_i]}$ for $t \in [\tau_i, \tau_{i+1}]$ remain to be solved over each subinterval, using equations \eqref{subeq:CMM} and \eqref{subeq:sourceAccum} respectively.

\newcommand{\CoarseGrid}{G}
\newcommand{\Xf}{G_\mathrm{f}}
\newcommand{\XPsi}{G_\psi}
\section{Numerical Implementation}
\label{sec:numimp}

The numerical implementation of the CMM without source term has been described in previous work \cite{yin2021characteristic,yin2021thesis, yin2023characteristic}. This section will present the necessary extensions to handle the source term.

We quickly recall some notations. We denote by $\vhX$ and $\hF$ the numerical discretizations of the (sub)map $\vX$ and (sub)source accumulation $F$. The spatial interpolation operators for these fields are $\Ih_M$ and $\Ih_A$ respectively, on grids $M$ and $A$. At time $t_n$, we denote by $\tilde{\vtx}_n$ the numerical vorticity evaluated from equation \eqref{subeq:vortEval}. The numerical velocity field, defined in space through an interpolant, is $\tilde{\vu}_n =  K_l \tilde{\vtx}_n$, where $K_l$ is a regularized numerical operator approximating the exact Biot-Savart operator $K$ up to some length scale $l$. The velocity field $\tilde{\vu}(\vx, t)$ for $t \in [t_n, t_{n+1})$, up to the next time step is defined using a temporal extrapolation operator $E_{\gamma} : (\tilde{\vu}_{n-\gamma+1} (\vx), \tilde{\vu}_{n-\gamma+2} (\vx), \ldots, \tilde{\vu}_{n}(\vx)) \mapsto \tilde{\vu}(\vx, t)$. The backward one-step map from this velocity is $\vtX_{[t_{n+1}, t_n]}(\vx)$ defined at $\vx$ pointwise by the solution of the backward in time ODE with velocity field $\tilde{\vu}$. Numerically, this is approximated by a classical Runge-Kutta-4 scheme.

In the absence of a source term, the CMM for incompressible Euler equations in 2d is summarized by an iteration of the following steps,
\begin{subequations}
    \begin{align}
        & \tilde{\vu}_n = K_l ( \vtx_0 \circ \vhX_{[t_n, \tau_i]} \circ  \vhX_{[\tau_i, \tau_{i-1}]} \circ \cdots \circ \vhX_{[\tau_1, 0]}) \\
        & \tilde{\vu} = E_{\gamma} (\tilde{\vu}_{n-\gamma+1} , \tilde{\vu}_{n-\gamma+2} , \ldots, \tilde{\vu}_{n}) \\
        & \partial_r \vtX_{[t_{n+1}, r]} = - \tilde{\vu}(\vtX_{[t_{n+1}, r]}, r) , \qquad \vtX_{[t_{n+1}, t_{n+1}]} (\vx)  =  \vx \\
        & \vhX_{[t_{n+1}, \tau_i]} = \Ih_M [ \vhX_{[t_{n}, \tau_i]}  \circ \vtX_{[t_{n+1}, t_n]} ].
    \end{align}
\end{subequations}

\subsection{Numerical Implementation of the Source Term}

In the case of equations with source terms, during a given time subinterval $[\tau_i, \tau_{i+1}]$, the accumulated source term $F_{[t, \tau_i]}$ must be computed by discretizing equation \eqref{eq:sourceTermAccumEq}. Alternatively, given the flow maps $\vX_{[t, s]}$ of the velocity field $\vu$, we may write $F_{[t_{n+1}, \tau_i]}$ explicitly using Duhamel's principle as a function of the map, $f$ and $F_{[t_n, \tau_i]}$:
\begin{equation}
    F_{[t_{n+1}, \tau_i]} = F_{[t_n, \tau_i]} \circ  \vX_{[t_{n+1}, t_n]}  + \int_{t_n}^{t_{n+1}} f ( \vX_{[t_{n+1}, s]}) , s) ds .
\end{equation}

This gives us a direct semi-Lagrangian discretization formula:
\begin{equation}
    \hF_{[t_{n+1}, \tau_i]} = \Ih_A \left[ \hF_{[t_n, \tau_i]} \circ  \vtX_{[t_{n+1}, t_n]}  + \dt \sum_{j=0}^{k-1} a_j \tilde{f} ( \vtX_{[t_{n+1}, s_j]}) , s_j)  \right] , \label{eq:sourceTerm_semiLagrangian}
\end{equation}
where $a_j, s_j$ are $k$-stage the Runge-Kutta time quadrature weights and loci. In our current implementation, we opted to use the same Runge-Kutta scheme for the one-step map and source term integrations. Therefore, the one-step map $\vtX_{[t_{n+1}, t_n]}$ and its intermediate stage values $\vtX_{[t_{n+1}, s_j]}$ were already computed during the map integration and therefore may be reused at no additional cost.

The source term $\tilde{f}$ is also defined in $[t_{n}, t_{n+1}]$ using the time extrapolation operator $E_{\gamma}$ from discrete time data $\tilde{f}_n$. At each time $t_n$, the source term, in this case, the curl of the Lorentz force, is approximated by sampling the magnetic field using expression \eqref{subeq:magneticField_pullback} on the grid $A$. Then, the magnetic field is filtered in Fourier space to retain spatial scales well-resolved by the grid, from which we then compute the curl of the Lorentz force, giving us the source term $\tilde{f}_n$. We note that this operator is solely a function of $\vB_0$ and the map $\vhX_{[t, 0]}$, hence we denote $\tilde{f} = S_l (\vhX_{[t, 0]})$ where $S_l$ is given by $ \curl \left(\curl  \tilde{\vB} \times \tilde{\vB} \right)$ where $\tilde{\vB}$ is a filtered version of $\adj(\nabla\vhX_{[t_n, 0]}) \vB_0 (\vhX_{[t_n, 0]})$.

Here the method is described using two separate grids $M$ and $A$ for the map and the source term accumulation respectively. In addition to accuracy considerations, these grids' resolution also controls the remapping frequency. Indeed, since the evolution of both $\vhX$ and $\hF$ are discretizations of an advection equation on a fixed Eulerian grid, convergence is lost when the numerical solution reaches the Nyquist frequency of the grid. Therefore a rule of thumb is to trigger a remapping routine whenever either $\vhX$ and $\hF$ develops significant small scales ill-resolved by their respective grids $M$ and $A$. The freedom of using different grids then allows for further optimization between the remapping frequencies of the two terms. Generally, using the same grid should be sufficient since both terms are advected by the same velocity and minor tweaking is possible depending on the respective right-hand-sides for the two equations.

The CMM for the ideal MHD equations in 2d is summarized in the pseudocode algorithm \ref{algo:jfb}. In this algorithm, we have introduced two distinct filtering scales $l_1$ and $l_2$ for the operators $K$ and $S$ respectively. The filtering $l_1$ is chosen to guarantee that the resulting velocity field $\tilde{\vu}$ is well resolved on the map grid. Indeed, any subgrid scales of the velocity in the map evolution will only appear as noise since it cannot be captured by the interpolant $\Ih_M$. This filtering allows us to access, for each fixed $l_1$, the convergence regime in the limit of fine spatial and temporal grids. Similarly, $l_2$ is chosen so that the source term is well resolved by $\Ih_A$. The convergence of the numerical solution to the exact solution then relies on the $l_1, l_2$ limit which one should take after the spatial and temporal grid limits. A diagonalization argument may then be used to select a single subsequence. We give an overview of the error estimate arguments in the following section.

\begin{algorithm}[htp!]
    \small
    \begin{algorithmic}[1] 
        \renewcommand{\algorithmicrequire}{\textbf{Input:}}
        \renewcommand{\algorithmicensure}{\textbf{Output:}
        $\vhX_{[t_{n+1}, \tau_i]},\hF_{[t_{n+1}, \tau_i]}$}
        \REQUIRE{Current submap $\vhX_{[t_n, \tau_{i-1}]}$ sub integral $\hF_{[t_n, \tau_{i-1}]}$,}
        \REQUIRE{
        previous velocities $\tilde{\vu}_{n-1},\dots,\tilde{\vu}_{n-\gamma}$  and source terms
        $\tilde{f}_{n-1} ,\dots, \tilde{f}_{n-\gamma}$
        }
        \REQUIRE{Stack of submaps $\{\vhX_{[\tau_{i}, \tau_{i-1}]}\}_{i}$ and sub-integrals $\{\hF_{[\tau_{i}, \tau_{i-1}]}\}_{i}$}
        \ENSURE{}
        \STATE{$i=I_\text{maps}$}
        \WHILE{$i > 1$}
        \STATE{$\vhX_{[t_n, \tau_{i-1}]} = \vhX_{[\tau_{i}, \tau_{i-1}]} \circ \vhX_{[t_n, \tau_{i}]}$} 
        \COMMENT{follow characteristics backward in time}
        \STATE{$\hF_{[t_n, \tau_{i-1}]} = \hF_{[t_n, \tau_{i}]} + \hF_{[\tau_{i}, \tau_{i-1}]} \circ \vhX_{[t_n, \tau_{i}]}$} 
        \COMMENT{Add up sub-integrals}
        \STATE{$i=i-1$}
        \ENDWHILE
        \STATE{$\tilde{\vu}_n = K_{l_1} ( \vtx_0 \circ \vhX_{[t_n, 0]} + \hF_{[t_n, 0]} )$} 
        \COMMENT{Solve Biot-Savart at current time instant $n$}
        \STATE{$\tilde{f}_n = S_{l_2} ( \vhX_{[t_n, 0]})$} 
        \COMMENT{Compute source term at current time instant $n$}
        \STATE{
         $\tilde{\vu} = E_{\gamma} (\tilde{\vu}_{n-\gamma+1} , \tilde{\vu}_{n-\gamma+2} , \ldots, \tilde{\vu}_{n}), \quad \tilde{f} = E_{\gamma} (\tilde{f}_{n-\gamma+1} , \tilde{f}_{n-\gamma+2} , \ldots, \tilde{f}_{n}),$
        }
        \COMMENT{Lagrangian Interpolants}
        \STATE{ $\vtX_{[t_{n+1}, t_n]}(\vx) = \vx - \dt  \sum_{j=0}^{k-1} a_j \tilde{\vu} ( \vtX_{[t_{n+1}(\vx), s_j]}) , s_j)$} 
        \COMMENT{Runge-Kutta stepping of map}
         \STATE{$\vhX_{[t_{n+1}, \tau_i]} = \Ih_M \left[ \vhX_{[t_{n}, \tau_i]}  \circ \vtX_{[t_{n+1}, t_n]} \right],$}
         \COMMENT{Compose submap}
        \STATE{$ \hF_{[t_{n+1}, \tau_i]} = \Ih_A \left[ \hF_{[t_n, \tau_i]} \circ  \vtX_{[t_{n+1}, t_n]}  + \dt \sum_{j=0}^{k-1} a_j \tilde{f} ( \vtX_{[t_{n+1}, s_j]}) , s_j)  \right] $}
        \COMMENT{Runge-Kutta Integration of source-term}
        \RETURN{}
    \end{algorithmic}
    \caption{\texttt{CMM\_timestep()} - Pseudo code of the CMM-MHD method.}
    \label{algo:jfb}
\end{algorithm}

\subsection{Error Analysis}
\label{subsec:ErrorAnalysis}

In this section we present a sketch of the error estimates for the CMM with source term, detailed error analysis for the CMM can be found in \cite{TaylorYinNave2024functional}. One notable difference is that here, instead of the momentum equation, we use the vorticity equation which, in the 2d incompressible case, is simply a scalar advection equation, thereby avoiding derivative loss due to the differential pullback operator. Under suitable assumptions on the spatial and temporal discretization operators, the following error estimate can be reached.

\begin{proposition}
    The CMM for the ideal MHD equations in two dimensions using Hermite cubic spatial interpolation and $3^{rd}$ order Lagrange time extrapolation and Runge-Kutta 4 integration has a numerical error of order $\mathcal{O}(N^{-3} + \dt^{3} ) + o(l_1 + l_2)$ where $h$ and $\dt$ are the spatial and temporal grid spacing and $l_1, l_2$ are the filtering length scales of the modified equation.
\end{proposition}

To simplify the analysis, we assume that the spatial discretization $\Ih_M$ for the map and $\Ih_A$ for the source term accumulation, on grids of cell width $N^{-1}$, are smooth projections on Sobolev spaces $H^s$ for some $s > n/2 +1$ with the contraction and approximation properties
\begin{equation}
    | \Ih [f ] |_r \leq | f |_r, \qquad \| f - \Ih [f] \|_r \leq C | f |_s N^{-(s-r)}
\end{equation}
for any $0 \leq r \leq s$, for instance, Fourier projections on $N^2$ harmonics \cite{boyd2001chebyshev} or mollification to scale $N^{-1}$. We also assume that once the numerical velocity $\vtu$ is defined in $[t_n, t_{n+1}]$ as a Lagrange polynomial, that the numerical integration for the one-step map $\vtX_{[t_{n+1}, t_n]}$ is performed exactly. In practice, the numerical solution obtained using Hermite cubic interpolants and RK4 integrators can then be viewed as a perturbation of the smooth scheme, with perturbation errors given by the accuracy of the interpolation and integration schemes. 

We define $\vthX_{[t, 0]} = \vhX_{[t_n, 0]} \circ \vtX_{[t, t_n]}$, thus satisfying
\begin{align}
    & ( \partial_t + \vtu \cdot \nabla ) \vthX_{[t, 0]} = 0, \quad t \in (t_n, t_{n+1}], \\
    & \vthX_{[t_n, 0]} = \vhX_{[t_n, 0]} .
\end{align}
We may then define a time-interpolated map $\vhX_{[t, 0]} = \Ih_M[ \vthX_{[t, 0]}]$ for all $t$, and this map corresponds to the numerical map obtained from the method at discrete times $t=t_n$. It follows from the linearity of the interpolation operator that there exists some velocity field $\vv_{\vhX} \in L^1(0, t, H^{s-1})$ which generates the map $\vhX_{[t, 0]}$. Since we may write $\vhX_{[t, 0]}$ exactly as the backward flow map of some velocity field, it follows that advected fields obtained by pullback through this map also satisfy an advection equation. Letting $\tilde{\theta}(t) = \theta_0 \circ \vhX_{[t, 0]}$, then $\tilde{\theta}$ solves
\begin{equation}
    (\partial_t + \vv_{\vhX} \cdot \nabla ) \tilde{\theta} = 0, \qquad \tilde{\theta} (0) = \theta_0 .
\end{equation}
Let $\theta(t)$ be the exact solution of the advection of the same initial field $\theta_0$ by the velocity field $\vu$. Applying Gr\"{o}nwall's lemma we have that for $r \leq s-1$,
\begin{equation}
    \| \tilde{\theta} (t) - \theta (t) \|_r \leq C \| \nabla \theta_0 \|_{s-1} \| \vv_{\vhX} - \vu \|_{L^1(0, t, H^r)} .
\end{equation}

In terms of the velocity error, we have that
\begin{equation}
    \| \vv_{\vhX} - \vu \|_r \leq \| \vu - \vtu \|_r + C N_M^{-(s-r-1)},
\end{equation}
where $N_M^{-1}$ is the cell width for the map grid $M$. The temporal error $\| \vu - \vtu \|_r$ contains the nonlinearity as the numerical velocities $\vtu$ depend on the numerical map. The numerical velocity $\vtu$ is defined through an order $\gamma$ extension operator $E_{\gamma}$ which extrapolates a velocity $\vtu \in L^1(t_n, t_{n+1}, H^s)$ from given numerical velocity data at steps $n-\gamma +1$ to $n$. We denote by $\vtu_i$ the numerical velocity at $t_i$ and $\vu(t_i)$ the exact velocity. The accuracy and boundedness properties
\begin{align}
    & \| \vv(t) - E_{\gamma} (\vv(t_{n-\gamma+1}), \ldots, \vv(t_n) \|_r \leq \| \partial_t^{\gamma+1} \vv \|_r \dt^{\gamma+1}, \\
    & \| E_{\gamma} (\vv_0, \ldots, \vv_{\gamma-1}  ) \|_r \leq C \max_{i=0, \ldots, \gamma-1} \| \vv_i \|_r ,
\end{align}
are sufficient to control the time approximation error:
\begin{align}
    \| \vu - \vtu \|_r & \leq \| \vu(t) - E_{\gamma} (\vu(t_{n-\gamma+1}), \ldots, \vu(t_n) \|_r \\
    & + \| E_{\gamma} ( (\vu(t_{n-\gamma+1}) - \vtu_{t_{n-\gamma+1}}) , \ldots,  (\vu(t_{n}) - \vtu_{t_{n}})  )  \|_r \\
    & \leq C \dt ( \dt^{\gamma} + \max_{i = 0, \ldots, \gamma-1} \| \vu(t_{n-i}) - \vtu_{t_{n-i}}) \|_r ) .
\end{align}
It remains to control the discrete step velocity error $\vu(t_i) - \vtu_{t_i}$ which involves the map pullback, Biot-Savart kernel, and the source term integration.

The numerical accumulated source term $\hF_{[t_n, 0]}$ is obtained from semi-Lagrangian methods using velocity $\vtu$ and source term $\tf$. We again try to write an advection equation for $\hF_{[t,0]} = \Ih_A [ \tilde{\hF}_{[t, 0]} ]$:
\begin{equation}
    (\partial_t + \vtu \cdot \nabla) \hF_{[t, 0]} = \tf - [ \Ih_A, \vtu \cdot \nabla ] \tilde{\hF}_{[t, 0]} ,
\end{equation}
where $[ \Ih_A, \vtu \cdot \nabla ]$ denotes the commutator between interpolation and the material derivative. This commutator arises from the fact that the source term is computed on a fixed grid rather than directly through the map and therefore subgrid scales generated by the advection are lost. The $L^1(0, t, H^q)$ error of this term is $\mathcal{O}(h^{s-q-1})$, i.e. the order of the scheme, but deteriorates accuracy once $\hF_{[t, 0]}$ contains significant small-scale features, at which point remapping is needed in order to reset this commutator error to 0.

Applying Gr\"{o}nwall yields the following estimate for the source term accumulation:
\begin{equation}
    \| \hF_{[t, 0]} - F_{[t, 0]} \|_q \leq \| \nabla \hF_{[t, 0]} \|_{s-1} \| \vtu - \vu \|_{L^1(0, t, H^q)} + C N_A^{-(s-q-1)} + \| f - \tf \|_{L^1(0, t, H^q)} .
\end{equation}
Again, defining $\tf(t)$ through the extension operator $E_{\gamma}$ gives an error 
\begin{equation}
    \| f - \tf \|_{L^1(0, t, H^q)} \leq C t ( \max_{i=1, \ldots, n} \| f(t_i) - \tf_i \|_q + \dt^\gamma ) .
\end{equation}
The control one can expect for the numerical source term depends heavily on the equation studied. For the ideal MHD equations and more generally for so-called equations with advected parameters, it turns out that $f$ is solely a function of the backward map $\vX_{[t, 0]}$ (if we consider initial data as parameters). Indeed, for the curl of the Lorentz force
\begin{equation}
     \curl \left(\curl  \vB \times \vB \right) =  \curl \left(\curl  (\adj(\nabla\vX_{[t, 0]}) \vB_0 (\vX_{[t, 0]}) ) \times \adj(\nabla\vX_{[t, 0]}) \vB_0 ( \vX_{[t, 0]}) \right) ,
\end{equation}
we may abstractly write a source term operator $S$ as $f(t) = S ( \vX_{[t, 0]})$. This can then be approximated by a regularized numerical operator $S_{l_2} : \textit{SDiff}^s \to H^q$ which should be $C^1$ so that $\| S_{l_2} (\vX_{[t, 0]}) - S_{l_2} ( \vhX_{[t, 0]} ) \|_q \leq C \| \vv_{\vhX} - \vu \|_{L^1(0, t, H^q)}$. The source term error can then be controlled by
\begin{equation}
    \| f(t_i) - \tf_i \|_q \leq \| (S_{l_2} - S) ( \vX_{[t_i, 0]}) \|_q + C \| \vv_{\vhX} - \vu \|_{L^1(0, t, H^q)},
\end{equation}
where the first term is the error of the approximate source term operator along the exact solution. Therefore,
\begin{equation}
    \| \hF_{[t, 0]} - F_{[t, 0]} \|_q \leq C( N_A^{-(s-q-1)} + \dt^{\gamma} +  \| (S_{l_2} - S) ( \vX_{[t_i, 0]}) \|_q + \| \vv_{\vhX} - \vu \|_{L^1(0, t, H^q)} ) .
\end{equation}

We now have everything necessary for estimating $ \| \vu(t_n) - \vtu_n \|_r$. Denote by $K$ the Biot-Savart kernel $- \nabla \times \Delta^{-1} : H^{r-1} \to H^{r}$ and let $K_l$ be a regularized numerical approximation $K_l : H^{r-k} \to H^r$ for $k \geq 1$. We have the exact and numerical velocities are given by
\begin{equation}
    \vu(t_n) = K \omega (t_n) = K (\omega_0 \circ \vX_{[t_n, 0]} + F_{[t_n, 0]}) , \qquad \vtu_n =  K_{l_1} \omega_n = K_l ( \omega_0 \circ \vhX_{[t_n, 0]} + \hF_{[t_n, 0]} ) .
\end{equation}
We obtain the following bound on the error:
\begin{align}
    & \| \vu (t_n) - \vtu_n \|_r  \lesssim \| ( K - K_{l_1} ) \omega (t_n) \|_r + \| \omega_0 \circ \vX_{[t_n, 0]} - \omega_0 \circ \vhX_{[t_n, 0]} \|_{r-k} + \| F_{[t_n, 0]} - \hF_{[t_n, 0]} \|_{r-k} \\
    & \leq \| ( K - K_{l_1} ) \omega (t_n) \|_r + \| (S-S_{l_2}) ( \vX_{[t_i, 0]}) \|_{r-k} \\
    & + N_M^{-(s-r-1)} + N_A^{-(s-r+k-1)}+ \dt^\gamma + \| \vtu - \vu \|_{L^1(0, t_n, H^{r})} ,
\end{align}
where we omitted any constant factors depending on the exact solution and the bounds on the $H^s$ norms of the numerical maps and the $H^{s-1}$ norms on the numerical source term accumulation $\hF$. Control on the growth of $\| \vhX_{[t, 0]} \|_s$ is studied in \cite{TaylorYinNave2024functional}, the growth rate of $\| \hF_{[t, 0]} \|_{s-1}$ relies on the choice of numerical solver for the advection equation. Both of these can be limited using the remapping method. In fact, since all constants in the error estimates are controlled by the $s$ seminorm of the map and $s-1$ seminorm of the accumulated source term, this suggests that the decay rate of the Fourier spectrum of the numerical map and source term data are good criteria for triggering the remapping step.

A last application of Gr\"{o}nwall's inequality then gives us an order $s-r-1$ in space and $\gamma$ in time convergence for the numerical velocities and by extension the map, vorticity and magnetic fields. We note that in the $(N^{-1}, \dt) \to 0$ asymptotic, the leading error is the operator errors $K-K_{l_1}$ and $S - S_{l_2}$ which depend on the filtering parameters ${l_1}, {l_2}$. This makes the error estimate highly non-trivial since the filtering parameters regularize the solution thereby limiting the presence of small scales in the numerical solution. On the one hand, choosing large ${l_1}, {l_2}$, i.e. a stronger regularization allows the error to enter the asymptotic regime for coarser grids, however, the deviation from the exact solution is also increased by the regularization errors $K-K_{l_1}$ and $S-S_{l_2}$. On the other hand, choosing smaller filtering scales reduces these errors but also amplifies the discretization errors due to potentially faster growth in the $H^s$ norms of the maps. The optimal relation between ${l_1}, {l_2}$ and $N^{-1}, \dt$ is highly dependent on the exact solution. As a general rule of thumb, for a given filter, $N^{-1}$ and $\dt$ should be chosen sufficiently small so that all scales expected in an ${l_1},{l_2}$-regularized solution are well resolved. The exact solution can then be approached by taking the ${l_1},{l_2} \to \infty$ limit with optimal $N^{-1}, \dt$ for each given filter. We note that the convergence of the numerical solution in the ${l_1}, {l_2} \to \infty$ limit might not be in a strong sense with respect to an $H^r$ norm. This remains a difficult problem and depends on whether the exact solution can be obtained as some weak limit of strong solutions to a family of regularized problems. The above shows that the regularization employed in this method is close to that of the Langrangian-averaged MHD equations or the related MHD-$\alpha$ equations for which global well-posedness is known \cite{linshiz2007analytical}. In this sense, the numerical results from the CMM, though fully inviscid, offer a relevant comparison with the Large Eddy Simulation models for which there is an extensive study in terms of spectral scaling and energy dissipation \cite{pietarila2006inertial, pietarila2009lagrangian, graham2005cancellation}.

In this paper, we employ a Hermite cubic spatial interpolation with third-order Lagrange extrapolation in time. This corresponds to a perturbation of a smooth scheme which is stable as long as the numerical solution $\vhX_{[t, 0]}$ is not highly oscillatory. More precisely, the error of this perturbation is controlled by the difference of $\vhX_{[t, 0]}$ with its Fourier projection $\Ih_{M,f} [\vhX_{[t, 0]} ]$ on the same grid $M$. This suggests that, in the case of Hermite cubic interpolation, we should take the $H^4$ norm of the Fourier transform of the map data as a remapping criterion. This holds similarly for $\hF_{[t, 0]}$.

\section{Numerical Studies} 
\label{sec:numtest}

For code validation of CMM with source terms, we first construct an exact solution of a linear advection equation with a swirling velocity field and a source term assessing the convergence properties. Then we apply CMM to the classical Orszag--Tang test case (OT) for ideal MHD and compare it with available results in the literature.

\subsection{Manufactured solution: linear advection with source term}
\label{subsec:manufactured_solution_linear_advection}

\begin{figure}[htp!]
    \centering
      \setlength\figureheight{0.85\textwidth}
     \setlength\figurewidth{0.95\textwidth}
    \includetikz{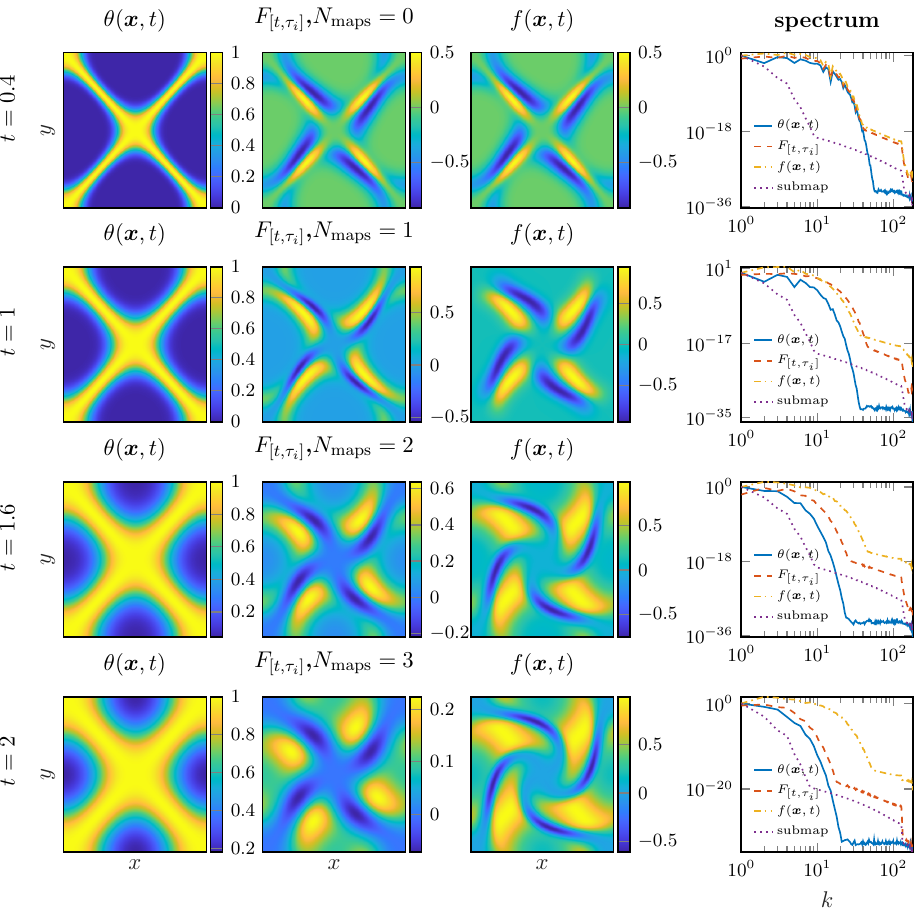}
    \caption{Linear advection swirl test case source term analysis. Shown is the solution $\theta(\vx, t)$, subintegrals $F_{[t,\tau_i]}$, source term $f(\vx, t)$ and corresponding Fourier spectra (from left to right), for $t=0.4, 1, 1.6$ and 2 (from top to bottom).}
    \label{fig:sourcetermFFT}
\end{figure}

In the first numerical example, we study linear advection with the source term
\begin{equation}
    \label{eq:advec_eq}
    \partial_t \theta + \vect{u}\cdot \nabla \theta = f
\end{equation}
using a divergence-free velocity field that generates a swirl:
\begin{align}
\vect{u}(x,y,t) =  \cos(t/4) \left(\sin^2(0.5x)\sin(y),-\sin(x)\sin^2(0.5y)\right)
\end{align}
on a periodic $[0,2\pi]^2$ domain and for $t\in[0,2]$. 
To test our framework, we impose a manufactured solution:
\begin{equation}
    \label{eq:manufactured}
    \theta^\text{ref}(x,y,t) = \exp\left( - \frac{( \cos(y) - \cos(x) )^2}{ (t-0.5)^2 + \epsilon } \right)
\end{equation}
with $\epsilon = 0.1 $ and where the corresponding source term $f$ is determined analytically by \cref{eq:advec_eq}. 
This test case is set up such that the solution has a constant spatial profile with narrowing crossing lines towards $t=0.5$ and increasing thickness after $t>0.5$. Therefore, the source term $f$ needs to compensate for the deformation created by the swirl velocity field.
The resulting evolution of the source term and solution is shown in \cref{fig:sourcetermFFT}. Additionally, we show the subintegrals $F_{[t,\tau_i]}$ for the different time instances and their Fourier spectra.
At time $t=0.4$ we see that the sub-integral and source are identical since no remapping was initiated.
Furthermore, the Fourier spectra of the source term/sub-integral and analytical solution $\theta$ are similar because of the small-scale structures present in the solution at time $t=0.4$.
In a source term free case, the scale in the solution would not affect the numerical precision because the evolved quantity is the map, not the state itself. 
Since in the presence of source terms, we have to store quantities that depend on the solution itself, the scales of the solution play a role in the numerical precision. Hence, to preserve these fine scales, storing the sub-integrals on finer grids might be necessary.
As soon as the first remapping is performed, we see that the spectra of the source term and sub-integrals differ. As the sub-integrals capture only local information of the solution, they have in general a faster decay than $f(x,t)$ itself. 
However, note that the spectrum of the submap decays for all time instances the fastest. 

\begin{figure}[htp!]
    \centering
      \setlength\figureheight{0.3\textwidth}
     \setlength\figurewidth{0.35\textwidth}
    \includetikz{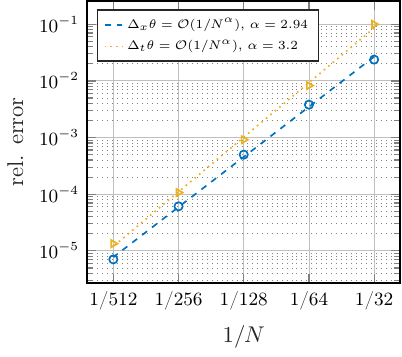}
    \caption{Convergence in space ($\Delta_x$) and time ($\Delta_t$) for the swirl test case, compared to the analytical reference solution.}
    \label{fig:convergence_space_time}
\end{figure}

Finally, we show the achieved third-order convergence in space and time in \cref{fig:convergence_space_time}.

Convergence in space and time is computed for $\Delta_{x/t}\theta =  \Vert\theta^\text{ref}(\cdot,T_\text{end})- \theta^N \Vert_\infty$ using the manufactured solution \cref{eq:manufactured}. Space convergence is achieved by increasing the number of grid points of the map $N=64,128,\dots,512$ in each dimension while keeping a constant time step of $\Delta t =1/512$ and integrating to $T_\text{end}=1$.
Convergence in time is achieved by decreasing the step size $\Delta t=1/64,1/128,\dots,1/512$, while keeping the number of grid points of the map at $N=512$ in each dimension. For both convergence tests, we keep the sampling grid of the velocity and the grid for evaluating all quantities at $512^2$.
For all convergence tests, remapping was switched off.
We remark that as mentioned before the fine-scale structure of the solution can not be decoupled from the numerical evolved quantities, as would be the case for the source-free swirl test case. Hence, we observe that for $\epsilon\ll 0.1$ the convergence rate is reduced, especially for the coarser resolutions.

\subsection{Orszag--Tang test case}

Our numerical experiments are based on the classical Orszag--Tang MHD test case that was implemented to study singularities in MHD turbulence \cite{orszag1979small}.

In this case, we expect $3^{\text{rd}}$ order convergence in the $L^2$ norm in space and time. In fact, our numerical experiments show $3^{\text{rd}}$ convergence in the $L^\infty$ norm. We define the following errors:
\begin{alignat}{3}
    \Delta \vmap &=  \Vert\vmap^\text{ref}(\cdot,t^n)- \vmap^n \Vert_{L^\infty}    &\quad& \text{(map error)}\,,                      \label{eq:map_error}\\
    \Delta \omega &=  \Vert \omega^\text{ref}(\cdot,t^n)- \omega^n \Vert_{L^\infty}               &\quad& \text{(vorticity function error) }\,,    \label{eq:dist_error}\\
    \Delta \Bfield &=  \Vert\Bfield^\text{ref}(\cdot,t^n)- \Bfield^n \Vert_{L^\infty}    &\quad& \text{(magnetic field error)}\,,                      \label{eq:map_error}
\end{alignat}
which will be studied numerically in \cref{sec:numtest}.

\begin{figure}
    \centering
     \setlength\figureheight{0.7\textwidth}
     \setlength\figurewidth{0.7\textwidth}
     \begin{subfigure}[b]{0.48\textwidth}
    \centering
      \setlength\figureheight{0.7\textwidth}
     \setlength\figurewidth{0.7\textwidth}
    \includetikz{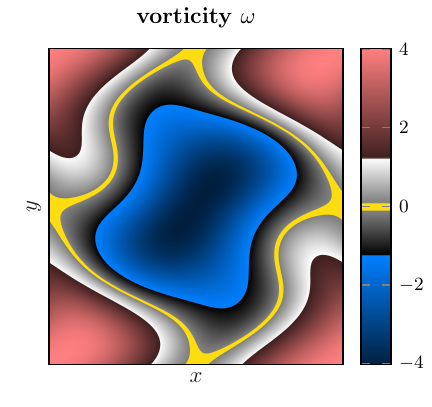}
    \end{subfigure}
    \begin{subfigure}[b]{0.48\textwidth}
    \centering
      \setlength\figureheight{0.7\textwidth}
     \setlength\figurewidth{0.7\textwidth}
    \includetikz{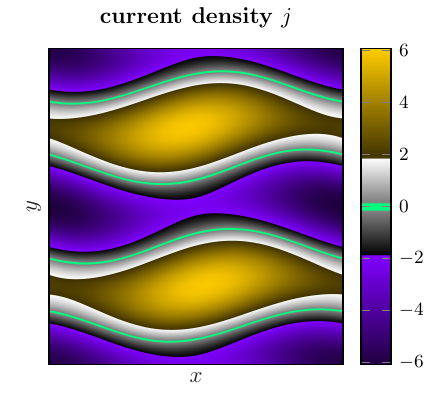}
    \end{subfigure}
    \begin{subfigure}[b]{0.48\textwidth}
    \centering
      \setlength\figureheight{0.7\textwidth}
     \setlength\figurewidth{0.7\textwidth}
    \includetikz{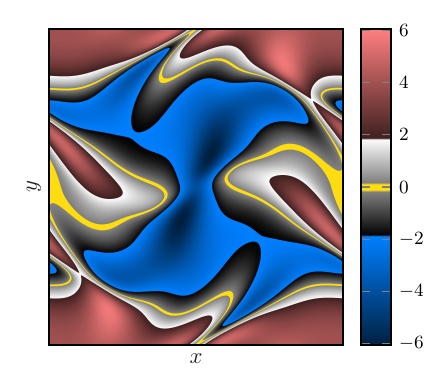}
    \end{subfigure}
    \begin{subfigure}[b]{0.48\textwidth}
    \centering
      \setlength\figureheight{0.7\textwidth}
     \setlength\figurewidth{0.7\textwidth}
    \includetikz{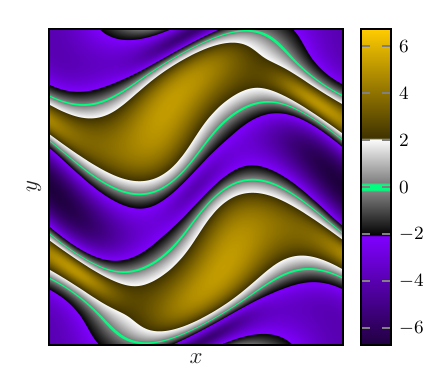}
    \end{subfigure}
    \begin{subfigure}[b]{0.48\textwidth}
    \centering
      \setlength\figureheight{0.7\textwidth}
     \setlength\figurewidth{0.7\textwidth}
    \includetikz{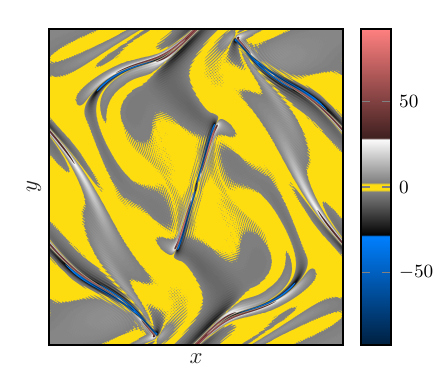}
    \end{subfigure}
    \begin{subfigure}[b]{0.48\textwidth}
    \centering
      \setlength\figureheight{0.7\textwidth}
     \setlength\figurewidth{0.7\textwidth}
    \includetikz{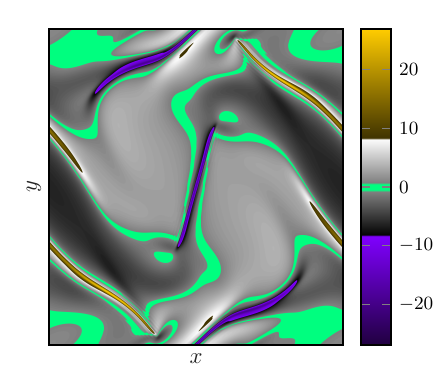}
    \end{subfigure}
    \caption{Time evolution of the vorticity $\omega$ (left) and current density $j$ (right) from top to bottom $t=0.1,0.4,0.9$ for OT.}
    \label{fig:MHD-time-evolve}
\end{figure}

\begin{figure}
    \centering
     \setlength\figureheight{0.7\textwidth}
     \setlength\figurewidth{0.7\textwidth}
     \begin{subfigure}[b]{0.48\textwidth}
    \centering
      \setlength\figureheight{0.6\textwidth}
     \setlength\figurewidth{0.7\textwidth}
    \includetikz{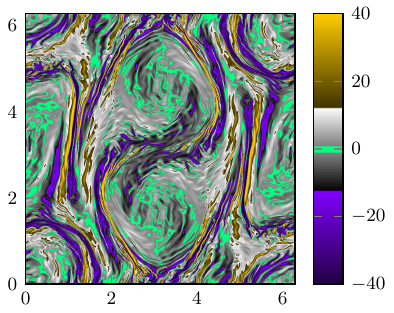}
    \end{subfigure}
    \begin{subfigure}[b]{0.48\textwidth}
    \centering
      \setlength\figureheight{0.6\textwidth}
     \setlength\figurewidth{0.7\textwidth}
    \includetikz{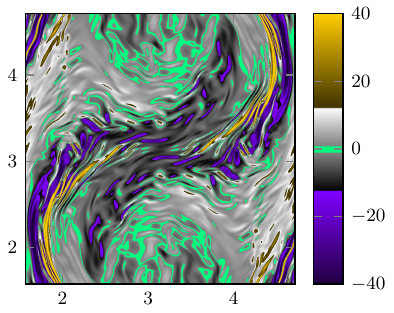}
    \end{subfigure}\\
    \begin{subfigure}[b]{0.48\textwidth}
    \centering
      \setlength\figureheight{0.6\textwidth}
     \setlength\figurewidth{0.7\textwidth}
    \includetikz{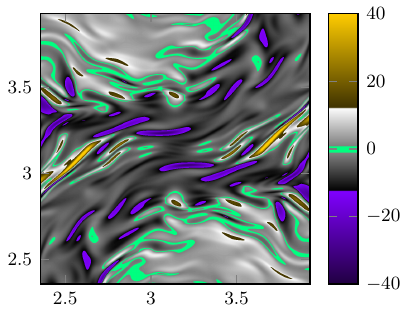}
    \end{subfigure}
    \begin{subfigure}[b]{0.48\textwidth}
    \centering
      \setlength\figureheight{0.6\textwidth}
     \setlength\figurewidth{0.7\textwidth}
    \includetikz{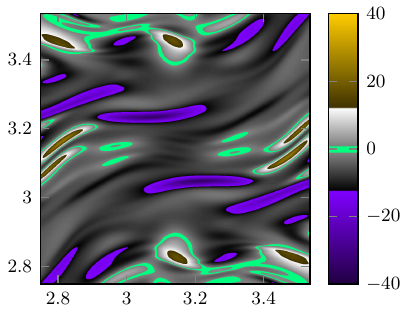}
    \end{subfigure}
    \caption{Zoom with center $(x,y)=(L/2,L/2)$ of the current density $j$ at time $t=4.0$ for OT. The zoom factor between each image is 2.}
    \label{fig:MHD-zoom-currentdisty}
\end{figure}

The initial condition of the test case reads:
\begin{align}
    \label{eq:inicond_OT}
    \omega(x,y,0) & = 2(\cos(2 x)+\cos(2 y))\,,\\
    \vect{B}(x,y,0) & = 2(-\sin(2 y), 2\sin(x))\,.
\end{align}
The test case is simulated inside the periodic domain $\Omega=[0,2\pi]^2$. A time evolution for $t = 0.1, 0.4, 0.9$ is shown for the vorticity $\omega$
and current density in \cref{fig:MHD-time-evolve}. Furthermore, successive zooms into the fine-scale structures of 
the current density are shown in \cref{fig:MHD-zoom-currentdisty} at time $t=4.0$.

For later reference, we define some characteristic quantities of the system.
The potential and kinetic energy is defined by:
\begin{equation}
    \Ekin = \frac12\norm{\vect{u}}^2_{L^2} \qquad \text{and} \qquad   \Epot = \frac12\norm{\vect{B}}_{L^2}\,.
\end{equation}
Here we use the standard $L_2$ norm $\norm{\vect{f}}^2=\langle \vect{f},\vect{f} \rangle$, with scalar product 
\begin{equation}
   \langle \vect{f},\vect{g}\rangle  = \frac{1}{(2\pi)^2} \int \vect{f}(\vect{x})\bar{\vect{g}}(\vect{x}) \d \vect{x}\,.
\end{equation}
The following quantities are conserved for sufficiently smooth solutions:
\begin{alignat}{3}
\Etot  &= \Ekin + \Epot  &\qquad &\text{(total energy)}\,,\label{Etot}\\
H_c  &=  \langle \vect{u},\vect{B}\rangle &\qquad &\text{(cross helicity)}\,,\label{eq:helicity}\\
A &= \frac{1}{2} \, \norm{\vect{A}}^2_{L^2} &\qquad &\text{(squared magnetic vector potential)}\,.\label{eq:squarePotential}
 \end{alignat}
where $\vect{B}=\mathrm{curl} \vect{A}$. The magnetic vector potential can be obtained from the current density $\vect{J} = \nabla \times \bm B$  using $\nabla^2 \vect{A} = \vect{J}$. Note, since all our computations are in 2d: $\vect{A}=a\vect{e}_z,\vect{J}=j\vect{e}_z$. 
Lastly, the magnetic and velocity spectra are given by:
\begin{align}
    \label{eq-def:fourier_spectra}
    E_{\vect{u}}(k) &= \frac{1}{2} \, \sum_{ k-1/2 < \norm{\bm k}_2 \le k+1/2} \norm{\widehat{\vect{u}}_{\bm k}}_2^2\\
    E_{\vect{B}}(k) &= \frac{1}{2} \, \sum_{ k-1/2 < \norm{\bm k}_2 \le k+1/2} \norm{\widehat{\vect{B}}_{\bm k}}_2^2\,
\end{align}
where $\widehat{\vect{u}}_k$ and $\widehat{\vect{B}}_k$ are the coefficients
of the Fourier transformed quantities $\vect{u}$ and $\vect{B}$, respectively.

\subsubsection{Convergence in space and time}
\label{subsec:Conv_in_space_and_time}
\begin{figure}
\begin{subfigure}[b]{0.48\textwidth}
    \centering
      \setlength\figureheight{0.7\textwidth}
     \setlength\figurewidth{0.7\textwidth}
    \caption{}
    \includetikz{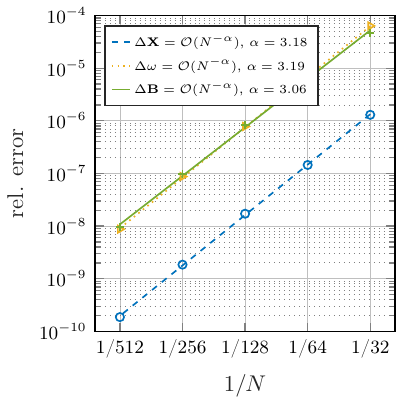}
    \end{subfigure}
    \begin{subfigure}[b]{0.48\textwidth}
    \centering
      \setlength\figureheight{0.7\textwidth}
     \setlength\figurewidth{0.7\textwidth}
     \caption{}
    \includetikz{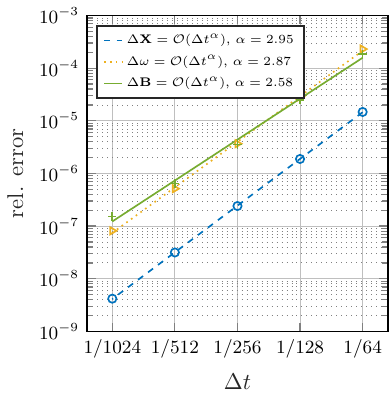}
    \end{subfigure}
    \caption{Convergence in space (a) and time (b) for OT. Shown is the relative error with respect to the reference solution in the $L^\infty$ norm.}
    \label{fig:conv_space}
\end{figure}

In \cref{subsec:manufactured_solution_linear_advection} we have studied the convergence in space and time for our numerical solution when compared to an analytic solution of the problem. Unfortunately, no analytic solution exists for ideal MHD. 
This is the reason why this section presents a self-comparison. Therefore our reference solution is computed with CMM using the finest resolution. For both, time and space convergence we use $N=2^{10}$ grid points in each dimension for the map and integrate up to time $t=0.1$ in time steps of size $\Delta t = 2^{-11}$. Furthermore, the resolution of the velocity field is kept at $1024\times 1024$, for all simulations.
To show numerical convergence in space we compare the reference solution to coarser approximations with lattice spacings $h\in 2\pi/2^{\{5,6,7,8,9\}}$. Similar for convergence in time we slowly increase the time steps from $\Delta t=2^{-6}$ to $\Delta t=2^{-10}$ and compute the $L^{\infty}$ error with respect to the reference solution. 
For both convergence tests remapping is switched off. Our experiments indicate 3rd-order convergence as shown in \cref{fig:conv_space} and \cref{table:EOC_OT}.

\begin{table}                                                         
\centering       
\begin{minipage}{.4\linewidth}
\begin{tabular}{cccc}                                                 
\toprule                                                              
$1/N$ & $\Delta \mathbf{X}$ & $\Delta \omega$ & $\Delta \mathbf{B}$ \\
\midrule                                                              
1/64 & 3.2 & 3.3 & 2.8 \\                                             
1/128 & 3.1 & 3.1 & 3.0 \\                                            
1/256 & 3.2 & 3.2 & 3.1 \\                                            
1/512 & 3.3 & 3.3 & 3.4 \\                                       
\bottomrule                                                           
\end{tabular}   
\subcaption{Spatial convergence}  
\end{minipage}%
\begin{minipage}{.4\linewidth}
\centering                                                                
\begin{tabular}{cccc}                                                     
\toprule                                                                  
$\Delta t$ & $\Delta \mathbf{X}$ & $\Delta \omega$ & $\Delta \mathbf{B}$ \\
\midrule                                                                   
1/128 & 3.0 & 3.0 & 2.9 \\                                                 
1/256 & 3.0 & 3.0 & 2.8 \\                                                 
1/512 & 2.9 & 2.8 & 2.5 \\                                                 
1/1024 & 2.9 & 2.7 & 2.1 \\                                             
\bottomrule                                                               
\end{tabular}                                                             
\subcaption{Temporal convergence}                                                  
\label{table:EOC_OT}                                                
\end{minipage}

\caption{Experimental order of convergence $
    \mathrm {EOC}(N)= {\log _2{\frac {\|f-f_{N}\|_{L^{\infty}}}{\|f-f_{N/2}\|_{L^{\infty}}}}}
    $ with reference solution computed at $N=2^{10}$ and $\Delta t = 2^{-11}$.}                                              
\label{table:MyTableLabel}                                            
\end{table}

\subsubsection{Comparison with literature results}
\label{subsec:comp_with_lit_results}
Next, we compare our results with a truncated Fourier Galerkin method applied to the ideal MHD equations at a resolution of $4096^2$ grid points~\cite{krstulovic2011alfven}. We use CFL time stepping with CFL number 1 and remapping threshold $\delta_\text{det} = 0.05$ for the comparison. The spatial resolution of the map is $512\times 512$, the velocity field as well as the field used for measurements is resolved at $1024\times 1024$. All parameters are summarized in \cref{tabl:reference_simulation}.

\begin{table}[htp!]
    \centering
    \begin{tabular}{l l}
    \toprule
        name                    & value \\
        \midrule
       Map resolution    & 512\\
       Velocity field resolution    & 1024\\
        inc. threshold $\delta_\text{det}$     & 0.05 \\
        time step size (CFL)                  & $1$\\
        cut off for map $l_1$ &         $(0.9k_\text{max})^{-1}$\\
        cut off for source term $l_2$ & $(0.1k_\text{max})^{-1}$\\
        \bottomrule
    \end{tabular}
    \caption{Parameters of the CMM for the reference computation.}
\label{tabl:reference_simulation}
\end{table}

\begin{figure}[htp!]
    \centering
      \setlength\figureheight{0.3\textwidth}
     \setlength\figurewidth{0.4\textwidth}
    \includetikz{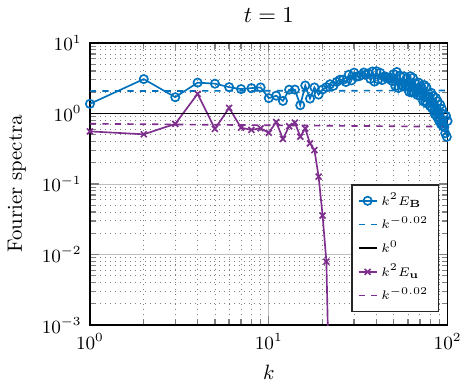}
    \caption{Compensated Fourier spectra $E_{\vect{u}}(k)$ and $E_{\vect{B}}(k)$ at $t=1$ and corresponding fit of the spectra using a wavenumber range from $k=1$ to $k=14$.}
    \label{fig:fourir_ub_1}
\end{figure}

For a first visual inspection, we compare the current density at time $t=0.9$ in the last row of \cref{fig:MHD-time-evolve} with Fig. 9a) in \cite{krstulovic2011alfven}.
We observe a similar structure. However, with closer inspection of the present simulations we observe the absence of some features present in \cite{krstulovic2011alfven}. Next, we compute the Fourier spectra $E_\vect{B},E_\vect{u}$ defined in \cref{eq-def:fourier_spectra} at $t=1$, which are shown in \cref{fig:fourir_ub_1}. We can confirm \cite{krstulovic2011alfven}, which observed decay of $E_\vect{B}\propto E_\vect{u}\propto k^{-2}$ in the initial dynamics. 

Additionally, we compare the time dynamics of the kinetic and potential energy, flux and Fourier mode ratios in \cref{fig:enter-label} to \cite{krstulovic2011alfven}. We see a good agreement of the dynamics up to $t=1$. Thereafter the two types of simulations diverge. 
\begin{figure}
\centering
\begin{subfigure}[b]{0.48\textwidth}
    \centering
      \setlength\figureheight{0.6\textwidth}
     \setlength\figurewidth{0.8\textwidth}
    \includetikz{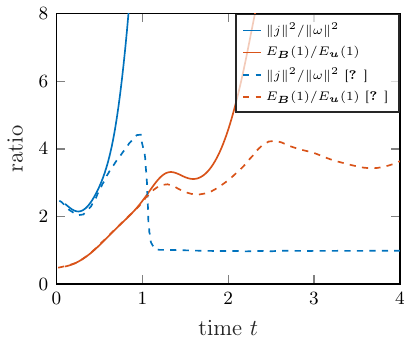}
    \caption{Flux and first Fourier mode ratios}
    \label{fig:enter-label}
\end{subfigure}%
\begin{subfigure}[b]{0.48\textwidth}
    \centering
      \setlength\figureheight{0.6\textwidth}
     \setlength\figurewidth{0.8\textwidth}
    \includetikz{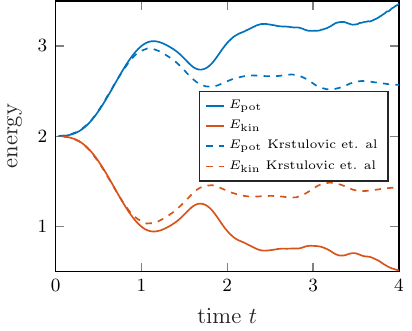}
    \caption{Potential and Kinetic Energy}
\label{fig:enter-label}
\end{subfigure}
\caption{Literature comparison to \cite{krstulovic2011alfven} using $N=512$ and $N_\psi=1024$.}
\label{fig:comparison-literature}
\end{figure}

\begin{figure}
    \centering
      \setlength\figureheight{0.6\textwidth}
     \setlength\figurewidth{0.8\textwidth}
    \includetikz{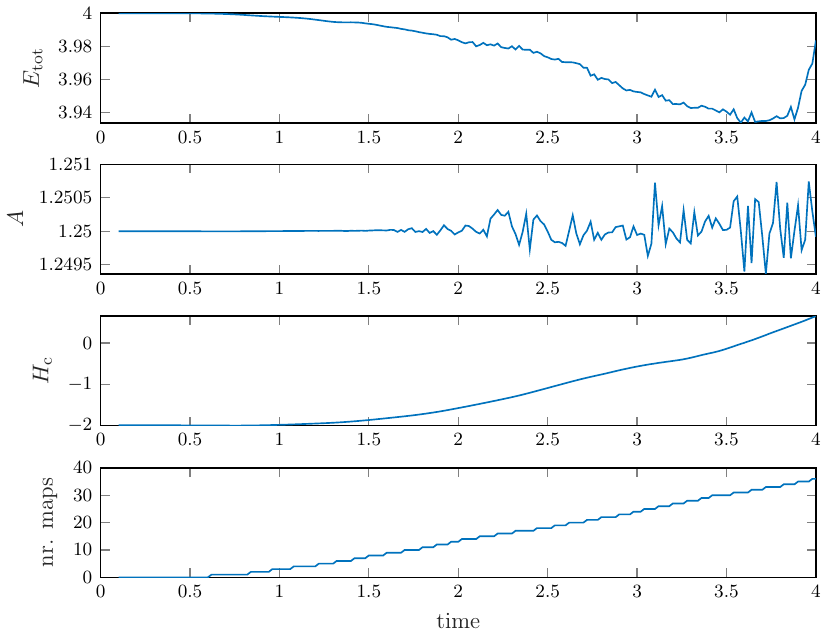}
    \caption{Conserved quantities of ideal MHD, total energy, squared magnetic vector potential and cross helicity (from top to bottom), and the number of remappings as a function of time for the OT test case computed at the finest resolution.}
    \label{fig:conserved_quantities}
\end{figure}

Lastly, \cref{fig:conserved_quantities} we plot the time evolution of the conserved quantities \cref{eq:helicity,Etot,eq:squarePotential} together with the number of submaps in \cref{fig:conserved_quantities}. Here we see a conservation of the quantities up to time $t=1$. Thereafter dissipation of total energy is observed. This might be caused by the formation of the thin current sheet which may imply a loss of regularity of the solution. 
Similar to our previous studies \cite{KrahYinBergmannSchneiderNave2023} we observe a linear growth in the number of submaps up to 40 maps at time $t=4$. 

\section{Conclusion}
\label{sec:concl}

In this present work, we presented an extension of the CMM for transport equations with source terms using the Duhamel integral. A recursive formula for the submap decomposition has been derived and its numerical implementation has been detailed. 
Numerical analysis showed global third-order convergence in space and time.
We validated the method in 2d and confirmed the theoretical third-order convergence in space and time using a manufactured solution. Then we benchmarked 2d ideal MHD, the classical Orzsag--Tang test case, and compared the results with truncated Fourier Galerkin computations. We likewise observed third-order convergence in space and time with respect to a fine grid reference computation. 

One unique feature of the CMM is its use of the structure of the diffeomorphism group to achieve an efficient decomposition of a long-time map using compositions of coarse grid submaps. Together with the preservation of the relabeling symmetry through the pullback operator, the method achieves arbitrary spatial resolution using coarse grid computations and the error is advective rather than diffusive in nature. These features are not straightforwardly preserved in the inhomogeneous case with source terms. In this work, this is remedied by applying a decomposition of the source term integration through Duhamel's formula. This extended methodology allows us to monitor the spatial resolution of the source term accumulation and subdivide the integral when numerical resolution is saturated, thereby controlling the numerical error arising from the transport term. These features are indeed observed in our numerical experiments. Using limited spatial resolution, the CMM was able to capture very small-scale features in the Orszag--Tang test case and we showed that subgrid scales are preserved by zooming into a subregion of the solution. This implies that the CMM is not prone to thermalization errors originating from a frequency-limited discretization of the quadratic convection term. In fact, our experiments could simulate longer periods without any apparent thermalization as observed in even higher-resolution pseudo-spectral codes.

Perspectives of the current work are the extension to MHD in three dimensions using a similar approach as done for the 3d incompressible Euler equations \cite{yin2023characteristic}.
Moreover, the developed source term technique will also enable us to consider kinetic equations, e.g. the Boltzmann equation with collision terms, which is a topic of current work \cite{KrahYinLinSchneiderNave2024}. 
Additionally, pursuing a possible coupling of MHD with kinetic kinetic equations, e.g. the developed Vlasov--Poisson CMM solver in \cite{KrahYinBergmannSchneiderNave2023} will provide efficient methods for solving the Vlasov--Maxwell system with possible application to magnetically confined fusion. 

Finally, on the application side and motivated by the singularity studies for 2d incompressible Euler equations using CMM~\cite{bergmann2024singularity}, we plan to extend these to incompressible MHD and seek numerical evidence whether solutions develop finite-time singularities, which is still an open question, see e.g.  \cite{cao2011global} where the global regularity has been analyzed theoretically in the presence of partial dissipation

\section*{Author Contribution Statement (CRediT)}

In the following, we declare the author's contributions to the publication:

\vspace{5pt}
{\small
\noindent
\begin{tabular}{@{}lp{10.8cm}}
\textbf{Xi-Yuan Yin:} &  implementation,  writing initial draft, numerical analysis, reviewing \& editing\\
\textbf{Philipp Krah:} & implementation, writing initial draft, numerical studies, visualization \\
\textbf{Jean-Christophe Nave:} & initial idea, editing and supervision\\
\textbf{Kai Schneider:} & initial idea, writing the initial draft, editing and supervision, funding acquisition, project administration \\
\end{tabular}
}




\section*{Acknowledgement}
The authors were granted access to the HPC resources of IDRIS under allocation No. AD012A01664R1 attributed by Grand Équipement National de Calcul Intensif (GENCI).
Centre de Calcul Intensif d’Aix-Marseille is acknowledged for granting access to its high-performance computing resources.
The authors acknowledge partial funding from the Agence Nationale de la Recherche (ANR), project CM2E, grant ANR-20-CE46-0010-01. JCN would like to acknowledge partial financial support from the NSERC Discovery Grant program.
The authors also thank Seth Taylor for many helpful discussions.


\bibliographystyle{abbrv}
\bibliography{references}

\end{document}